\documentclass[prx,amsfonts,amsmath,floatfix]{revtex4}

\usepackage{graphicx}
\usepackage{epsfig}
\usepackage{epstopdf}
\usepackage[usenames,dvipsnames]{xcolor}
\usepackage{braket}
\usepackage{hyperref}
\usepackage{amsthm}
\usepackage{enumitem}

\newcommand{\be}{\begin{equation}}
\newcommand{\ee}{\end{equation}}
\newcommand{\ba}{\begin{eqnarray}}
\newcommand{\ea}{\end{eqnarray}}

\newcommand{\bc}{\begin{center}}
\newcommand{\ec}{\end{center}}
\newcommand{\beq}{\begin{equation}}
\newcommand{\eeq}{\end{equation}}
\newcommand{\beqq}{\begin{equation*}}
\newcommand{\eeqq}{\end{equation*}}
\newcommand{\beqa}{\begin{align}}
\newcommand{\eeqa}{\end{align}}
\newcommand{\barr}{\begin{array}}
\newcommand{\earr}{\end{array}}
\newcommand{\bi}{\begin{itemize}}
\newcommand{\ei}{\end{itemize}}

\newtheorem{theo}{Theorem}

\begin{document}

\title{On the permanent of Sylvester-Hadamard matrices}

\author{Ulysse Chabaud$^{1}$}
\email{ulysse.chabaud@gmail.com}
\address{$^1$ Sorbonne Universit\'e, CNRS, Laboratoire d'Informatique de Paris 6, LIP6, 4 place Jussieu, 75115 Paris}
\date{\today}


\begin{abstract}
We prove a conjecture due to Wanless~\cite{wanless2005permanents} about the permanent of Hadamard matrices in the particular case of Sylvester-Hadamard matrices. Namely we show that for all $n\ge2$, the dyadic valuation of the permanent of the Sylvester-Hadamard matrix of order $n$ is equal to the dyadic valuation of $n!$. As a consequence, the permanent of the Sylvester-Hadamard matrix of order $n$ doesn't vanish for $n\ge2$.
\end{abstract}


\maketitle


Let $n\in\mathbb{N}^*$, the Sylvester-Hadamard matrix of order $n$ is defined by induction:
\begin{equation}
H_{n}=\begin{pmatrix} H_{n-1} & H_{n-1} \\ H_{n-1} & -H_{n-1} \end{pmatrix},
\end{equation}
with $H_0=1$. Equivalently,
\begin{equation}
H_{n}=\underbrace{H\otimes\dots\otimes H}_{n\text{ }times},
\end{equation}
where
\begin{equation}
H=H_1=\begin{pmatrix} 1 & 1 \\ 1 & -1 \end{pmatrix}.
\end{equation}
The Sylvester-Hadamard matrix of order $n$ is thus a square $(1,-1)$-matrix of size $2^n$. It is also a Hadamard matrix since $H_nH_n^T=2^nI_{2^n}$, where $I_{2^n}$ denotes the identity matrix of order $2^n$. While the entries of its first line and column are only $1$'s, each other line or column of $H_n$ contains the same number $2^{n-1}$ of $1$'s and $-1$'s. As a preliminary result, this implies that the products of the elements of each line or column of $H_n$ are equal to $1$ for $n\ge2$. Let $M\in\mathbb{N}^*$, let us define the permanent of an $M\times M$ matrix $T=(t_{ij})_{1\le i,j\le M}$ by
\begin{equation}
\text{Per}(T)=\sum_{\sigma\in \mathcal{S}_M}{\prod_{k=1}^{M}{t_{k\sigma(k)}}},
\end{equation}
where $\mathcal{S}_M$ is the symmetric group over $\{1,\dots,M\}$. For all $k\in\mathbb{Z}$ we will write $\nu_2(k)$ for the dyadic valuation of $k$. We may then state the following result:\\

\begin{theo} For all $n\ge2$, $\nu_2(\text{Per}(H_n))=2^n-1$.\\
\label{th:Sylv}
\end{theo}

\textbf{Proof.} For a given $n\ge2$, we write $H_n=S=(s_{ij})_{1\le i,j\le2^n}$. For $1\le i,j\le 2^n$, we denote by $S_{i,j}$ the matrix obtained from $S$ by removing the $i^{th}$ line and the $j^{th}$ column. We have
\begin{equation}
S=\underbrace {H\otimes\dots\otimes H }_{ n\text{ }times }.
\end{equation}
The lines of $H$ together with the element-wise multiplication form a group isomorphic to $\mathbb{Z}/2\mathbb{Z}$, thus the lines of $S$ together with the element-wise multiplication form a group isomorphic to $\left(\mathbb{Z}/2\mathbb{Z}\right)^n$. As a consequence, multiplying element-wise all lines of $S$ by its $k^{th}$ line, for any given $k\in\{1,\dots,2^n\}$, amounts to permuting the lines of $S$. For all $k\in\{1,\dots,2^n\}$, we can thus obtain the matrix $S_{1,1}$ from the matrix $S_{k,1}$ by multiplying element-wise all lines by the $k^{th}$ line of $S$ and then permuting the lines. Since the permanent is invariant by line permutation we obtain, for all $k\in\{1,\dots,2^n\}$,
\begin{equation}
\text{Per}(S_{k,1})=\epsilon_{k}\text{Per}(S_{1,1}),
\end{equation}
where $\epsilon_{k}=s_{1k}\prod_{l=1}^{2^n}{s_{kl}}$. For all $k\in\{1,\dots,2^n\}$, $s_{1k}=s_{k1}=1$, and the product of the elements of a line of $S$ is equal $1$ for each line with the preliminary result, so $\epsilon_{k}=1$. Hence for all $k\in\{1,\dots,2^n\}$,
\begin{equation}
\text{Per}(S_{k,1})=\text{Per}(S_{1,1}).
\label{pmper}
\end{equation}
We then use the Laplace expansion of the permanent to obtain
\begin{equation}
\begin{aligned}
\text{Per}(S)&=\sum_{k=1}^{2^n}{s_{k1}\text{Per}(S_{k,1})}\\
&=\sum_{k=1}^{2^n}{\text{Per}(S_{1,1})}\\
&=2^n\text{Per}(S_{1,1}),
\end{aligned}
\label{Laplace}
\end{equation}
where we used Eq.~(\ref{pmper}) in the second line. Now $S_{1,1}$ is a square matrix of size $2^n-1$, so by a result of~\cite{simion1983+}, $\nu_2(\text{Per}(S_{1,1}))=2^n-n-1$. For completeness we reproduce the proof here:
\begin{quotation}
Let $m\in\mathbb{N}^*$, let $A$ be any square $(1,-1)$-matrix of size $m$ and let $J$ be the square matrix of size $m$ whose entries are only $+1$. We write $A=J-2B$, for some $(0,1)$-matrix $B$. The expansion formula for the permanent of the sum of matrices~\cite{minc1984permanents} gives
\begin{equation}
\begin{aligned}
\text{Per}(A)&=\text{Per}(J-2B) \\
&=m!-2(m-1)!p_1(B)+2^2(m-2)!p_2(B)+\dots+(-1)^m2^mp_m(B),
\end{aligned}
\label{RHS}
\end{equation}
where $p_k(B)$ denotes the sum of the permanents of all the $k\times k$ submatrices of $B$, for $1\le k\le m$. Now for all $k$, $\nu_2(k!)=k-s_k$, where $s_k$ is the number of non-zero digits in the binary writing of $k$~\cite{koblitz2012p}. Hence for all $k<m$ we have 
\begin{equation}
\nu_2(2^{m-k}k!p_{m-k}(B))\ge (m-k)+(m-s_k) = m-s_k
\end{equation}
and $\nu_2(m!)=m-s_m$. If $m=2^n-1$, then $s_m=n$ while for all $k<m$ we have $s_k<n$, so $m-s_k>m-s_m$. This implies that the right-hand side of Eq.~(\ref{RHS}) has dyadic valuation $m-n=2^n-n-1$.
\end{quotation}

With Eq.~(\ref{Laplace}), this finally implies that
\begin{equation}
\nu_2(\text{Per}(S))=2^n-1.
\end{equation}
\qed

In particular, $\nu_2(2^n!)=2^n-1$, hence $2^{\nu_2(2^n!)}$ exactly divides $\text{Per}(H_n)$ for $n\ge2$, so the dyadic valuation of the permanent of the Sylvester-Hadamard matrix of order $n$ is equal to the dyadic valuation of $n!$ for $n\ge2$. A direct consequence of Theorem~\ref{th:Sylv} is that the permanent of the Sylvester-Hadamard matrix of order $n$ doesn't vanish for $n\ge2$. The latter result has been derived independently in~\cite{crespi2015suppression}.


\bibliographystyle{apsrev}
\bibliography{bibliography}

\end{document}